\documentstyle[11pt]{article}

\def\b#1{{\bf #1}}
\def\i#1{{\it #1}}

\begin{document}
\title{On the roots of the cubic defining the Tribonacci sequences}
\author{Mario Catalani\\
Department of Economics, University of Torino\\ Via Po 53, 10124 Torino, Italy\\
mario.catalani@unito.it}
\date{}
\maketitle
\begin{abstract}
\small{We consider a sequence of sums of powers of the the roots of the cubic
equation
characterizing the Tribonacci sequences and derive its relationship
with a particular Tribonacci sequence. Then we make a conjecture on the
possible generalization.}
\end{abstract}
\bigskip
\bigskip

\noindent
Let $\alpha,\,\beta,\,\gamma$ be the roots of the cubic $x^3-x^2-x-1=0$
(see for example \cite{mario3}, \cite{elia}, \cite{wolfram}).
Define
$$U_n=\sum_{i=1}^n\sum_{j=0}^{n-i}\alpha^i\beta^j\gamma^{n-i-j}.$$
We can write
$$U_{n+1}=\gamma U_n+\sum_{i=0}^n\alpha^i\beta^{n+1-i}+\alpha^{n+1}.$$
It follows
$$U_{n+2}=\gamma U_{n+1}+\sum_{i=0}^{n+1}\alpha^i\beta^{n+2-i}+\alpha^{n+2},$$
$$U_{n+3}=\gamma U_{n+2}+\sum_{i=0}^{n+2}\alpha^i\beta^{n+3-i}+\alpha^{n+3},$$
and
$$U_{n+4}=\gamma U_{n+3}+\sum_{i=0}^{n+3}\alpha^i\beta^{n+4-i}+\alpha^{n+4}.$$
We want to prove
$$U_{n+4}=U_{n+3}+U_{n+2}+U_{n+1}.$$
We are going to use mathematical induction. First of all $U_0=1$,
$$U_1=\alpha+\beta+\gamma=1,$$
$$U_2=\alpha^2+\beta^2+\gamma^2+\alpha\beta+\alpha\gamma+\beta\gamma=2,$$
$$U_3=\alpha^3+\beta^3+\gamma^3+\alpha^2\beta+\alpha\beta^2+\alpha^2\gamma
+\alpha\gamma^2+\beta^2\gamma+\beta\gamma^2+\alpha\beta\gamma=4.$$
These relationships follow from the fact $\alpha\beta\gamma =1,\,
\alpha\beta+\alpha\gamma+\beta\gamma=-1,\,\alpha+\beta+\gamma=1$. Then
$$U_3=U_2+U_1+U_0.$$
Now assume that
$$U_{n+3}=U_{n+2}+U_{n+1}+U_n.$$
Because $\alpha$ and $\beta$ are roots of the cubic $x^3-x^2-x-1=0$ we have
$\alpha^k= \alpha^{k-1}+\alpha^{k-2}+\alpha^{k-3}$,
$\beta^k= \beta^{k-1}+\beta^{k-2}+\beta^{k-3}$, $k=3,\,4,\,\ldots$. Then
$$U_{n+4}=\gamma (U_{n+2}+U_{n+1}+U_n)+
\sum_{i=0}^{n+3}\alpha^i\beta^{n+4-i}+
(\alpha^{n+3}+\alpha^{n+2}+\alpha^{n+1}).$$
Furthermore
\begin{eqnarray*}
\sum_{i=0}^{n+3}\alpha^i\beta^{n+4-i}&=&
\beta^{n+4}+\alpha\beta^{n+3}+\alpha^2\beta^{n+2}+\cdots
+\alpha^{n+3}\beta\\
&=&(\beta^{n+1}+\beta^{n+2}+\beta^{n+3})+\alpha(
\beta^n+\beta^{n+1}+\beta^{n+2})+\cdots\\
&&\quad\quad
+(\alpha^n+
\alpha^{n+1}+\alpha^{n+2})\beta\\
&=&\sum_{i=0}^n\alpha^i\beta^{n+1-i}+\sum_{i=0}^{n+1}\alpha^i\beta^{n+2-i}+
\sum_{i=0}^{n+2}\alpha^i\beta^{n+3-i}.
\end{eqnarray*}
Putting all together our claim is proved. So
$$U_n=U_{n-1}+U_{n-2}+U_{n-3},$$
with initial conditions $U_0=1,\,U_1=1,\,U_2=2$, as it can be easily seen.
Consider the Tribonacci sequence
$$T_n=T_{n-1}+T_{n-2}+T_{n-3},\qquad T_0=0,\,T_1=1,\,T_2=1.$$
It follows immediately
$$U_n=T_{n+1}.$$
In an analogous way it can be proved that if $\alpha$ and $\beta$ are the
roots of the quadratic $x^2-x-1=0$ characterizing Fibonacci sequences
($\alpha$ is the golden ratio) then
$$\sum_{i=0}^n\alpha^i\beta^{n-i}=F_{n+1},$$
where $F_n$ are Fibonacci numbers
$$F_n=F_{n-1}+F_{n-2},\qquad F_0=0,\,F_1=1.$$
Conjecture.

\noindent
Let $\alpha_i,\,i=1,\ldots,m$ be the roots of
$$x^m-x^{m-1}-x^{m-2}-\cdots -x-1=0.$$
Then
$$\sum_{i_1=0}^n \sum_{i_2=0}^{n-i_1}\cdots \sum_{i_{m-1}=0}^{n-i_1-\cdots
-i_{m-2}}\alpha_1^{i_1}\alpha_2^{i_2}\times\cdots\times
\alpha_{m-1}^{i_{m-1}}\alpha_m^{n-i_1-\cdots -i_{m-1}}
=V_{n+1},$$
where
$$V_n=V_{n-1}+V_{n-2}+\cdots +V_{n-m+1}+V_{n-m},$$
with initial conditions
$$V_0=0,\,V_1=V_2=\cdots =V_{m-2}=V_{m-1}=1.$$

\end{document}